\documentclass{article}
\usepackage{amsfonts,amsmath,amssymb}
\usepackage{graphicx,epsfig,psfig}
\usepackage{longtable,geometry}
\begin{document}


\newtheorem{theorem}{Theorem}[section]
\newtheorem{proposition}{Proposition}[section]
\newtheorem{lemma}{Lemma}[section]
\newtheorem{corollary}{Corollary}[section]
\newtheorem{remark}{Remark}[section]
\newtheorem{proof}{Proof:}

\renewcommand{\thesection}{\arabic{section}}
\renewcommand{\theequation}{\thesection.\arabic{equation}}
\renewcommand{\thetheorem}{\thesection.\arabic{theorem}}
\numberwithin{equation}{section}
\numberwithin{theorem}{section}
\numberwithin{proposition}{section}
\numberwithin{lemma}{section}
\numberwithin{remark}{section}
\setcounter{secnumdepth}{5}


\newcommand{\cl}{\centerline}
\newcommand{\sms}{\smallskip}
\newcommand{\ms}{\medskip}
\newcommand{\bs}{\bigskip}
\newcommand{\noi}{\noindent}
\newcommand{\itl}[1]{\textit{#1}}
\newcommand{\blf}[1]{\textbf{#1}}
\newcommand{\dsty}{\displaystyle}
\newcommand{\txty}{\textstyle}
\newcommand{\ssty}{\scriptstyle}
\newcommand{\tty}{\texttt}


\newcommand\Par{\mathhexbox278\,}


\newcommand{\al}{\alpha}
\newcommand{\Al}{\Alpha}
\newcommand{\be}{\beta}
\newcommand{\Be}{\Beta}
\newcommand{\Gm}{\Gamma}
\newcommand{\gm}{\gamma}
\newcommand{\dl}{\delta}
\newcommand{\Dl}{\Delta}
\newcommand{\lm}{\lambda}
\newcommand{\Lm}{\Lambda}
\newcommand{\kp}{\kappa}
\newcommand{\varep}{\varepsilon}
\newcommand{\vp}{\varphi}
\newcommand{\sig}{\sigma}
\newcommand{\Sig}{\Sigma}
\newcommand{\om}{\omega}
\newcommand{\Om}{\Omega}
\newcommand{\uom}{\mbox{\boldmath$\omega$}}
\newcommand{\btau}{\mbox{\boldmath$\tau$}}
\newcommand{\bnu}{\mbox{\boldmath$\nu$}}
\newcommand{\up}{\upsilon}
\newcommand{\z}{\zeta}


\newcommand{\df}[1]{\buildrel\mbox{\small def}\over{#1}}
\newcommand{\op}[1]{\buildrel\mbox{\tiny o}\over{#1}}
\newcommand{\db}{\prime\prime}
\newcommand{\bsl}{\backslash}
\newcommand{\lb}{\lbrack\!\lbrack}
\newcommand{\rb}{\rbrack\!\rbrack}
\newcommand\la{\langle}
\newcommand\ra{\rangle}
\newcommand{\ev}{\equiv}
\newcommand{\nev}{\not\equiv}
\newcommand{\nn}{\mathbb{N}}
\newcommand{\qq}{\mathbb{Q}}
\newcommand{\zz}{\mathbb{Z}}
\newcommand{\rr}{\mathbb{R}}
\newcommand{\rn}{\rr^N}
\newcommand{\cc}{\mathbb{C}}
\newcommand{\id}{\mathbb{I}}
\newcommand{\bo}{\mathbb{O}}

\newcommand{\amsb}[1]{\mathbb{#1}}
\newcommand{\mcl}[1]{\mathcal{#1}}
\newcommand{\bl}[1]{\mathbf{#1}}
\newcommand{\ov}[1]{\overline{#1}}
\newcommand{\wt}[1]{\widetilde{#1}}
\newcommand{\wh}[1]{\widehat{#1}}

\newcommand{\lra}{\longrightarrow}
\newcommand{\LLR}{\Longleftrightarrow}
\newcommand{\LRA}{\Longrightarrow}
\newcommand{\LLA}{\Longleftarrow}


\newcommand{\bbox}{\vrule height.6em width.6em 
depth0em} 
\newcommand{\os}{\vbox{\hrule \hbox{\vrule 
height.6em depth0pt 
\hskip.6em \vrule height.6em depth0em}
\hrule}} 


\newcommand{\dvg}{\operatorname{div}}
\newcommand{\curl}{\operatorname{curl}}
\newcommand{\supp}{\operatorname{supp}}
\newcommand{\essup}{\operatornamewithlimits{ess\,sup}}
\newcommand{\essinf}{\operatornamewithlimits{ess\,inf}}
\newcommand{\essosc}{\operatornamewithlimits{ess\,osc}}
\newcommand{\osc}{\operatornamewithlimits{osc}}
\newcommand{\sign}{\operatorname{sign}}
\newcommand{\loc}{\operatorname{loc}}
\newcommand{\diam}{\operatorname{diam}}
\newcommand{\dist}{\operatorname{dist}}
\newcommand{\card}{\operatorname{card}}
\newcommand{\meas}{\operatorname{meas}}
\newcommand{\spn}{\operatorname{span}}
\newcommand{\dtm}{\operatorname{det}}
%


\newcommand{\overlim}{\mathop{\overline{\lim}}\limits}
\newcommand{\underlim}{\mathop{\underline{\lim}}\limits}
\newcommand{\ttop}[2]{\genfrac{}{}{0pt}{}{#1}{#2}}
\newcommand{\bcu}{\mathop{\txty{\bigcup}}\limits}
\newcommand{\bca}{\mathop{\txty{\bigcap}}\limits}
\newcommand{\bsu}{\mathop{\txty{\sum}}\limits}
\newcommand{\pro}{\mathop{\txty{\prod}}\limits}


\newcommand{\pl}{\partial}
\newcommand{\ptt}{\frac{\pl}{\pl t}}
\newcommand{\ppx}{\frac\pl{\pl x}}
\newcommand{\dds}{\frac d{ds}}
\newcommand{\ddt}{\frac d{dt}}


\newcommand{\intl}{\int\limits}
\newcommand{\iintl}{\iint\limits}


\def\Xint#1{\mathchoice
    {\XXint\displaystyle\textstyle{#1}}%
    {\XXint\textstyle\scriptstyle{#1}}%
    {\XXint\scriptstyle\scriptscriptstyle{#1}}%
    {\XXint\scriptscriptstyle\scriptscriptstyle{#1}}%
    \!\int}
\def\XXint#1#2#3{\setbox0=\hbox{$#1{#2#3}{\int}$}
    \vcenter{\hbox{$#2#3$}}\kern-0.5\wd0}
\def\bint{\Xint-}
\def\dashint{\Xint{\raise4pt\hbox to7pt{\hrulefill}}}


\newcommand{\ovl}[3]{\int_{#1}^{#2}\kern-#3pt\raise4pt\hbox to7pt{\hrulefill}\ }

\newcommand{\ovll}[3]{\intl_{#1}^{#2}\kern-#3pt\raise4pt\hbox to7pt{\hrulefill}\ }

\newcommand{\tvl}[2]{\iint_{#1}\kern-#2pt\raise4pt\hbox to7pt{\hrulefill}\ }



\newcommand{\omt}{\Om_T}
\newcommand{\plo}{\partial\Omega}
\newcommand{\ovo}{\bar{\Om} }

%
\newcommand{\ci}[1]{C^\infty\!\left({#1}\right)}
\newcommand{\cio}[1]{C_o^\infty\!\left({#1}\right)}
\newcommand{\lloc}[1]{L_{\loc}\!\left({#1}\right)}
\newcommand{\xy}{|x-y|}


\newcommand{\intom}{\intl_{\Om}}
\newcommand{\intbo}{\intl_{\plo}}
\newcommand{\inom}{\int_{\Om}}
\newcommand{\inbo}{\int_{\plo}}
\newcommand{\intrn}{\intl_{\rn}}


\newcommand{\bye}{\end{document}}



%
%
%

%
%
%
%
%

%
%

\input harnack.mac
\begin{center}  {\huge\textbf{Regularity of $\infty$ for Elliptic Equations with Measurable
Coefficients and Its Consequences }}
\par \medskip\bigskip\end{center}
\begin{center} {\Large\textsc{Ugur G. Abdulla}}
\par \medskip\bigskip\end{center}
\begin{center} {\large\noindent \textsc{Department of Mathematics, Florida Institute of Technology, Melbourne, Florida 32901}}
\par \medskip\bigskip\end{center}
{\bf Abstract.} This paper introduces a notion of regularity (or irregularity)
of the point at infinity ($\infty$) for the unbounded open set $\Omega\subset\rr^{N}$ concerning second order uniformly elliptic equations with bounded and measurable coefficients, according as whether the  ${\cal A}$- harmonic measure of $\infty$ is zero (or positive). A necessary and sufficient condition for the existence of a unique bounded solution
to the Dirichlet problem in an arbitrary open set of $\rr^{N}, N\ge 3$ is established
in terms of the Wiener test for the regularity of $\infty$. It coincides with the Wiener 
test for the regularity of $\infty$ in the case of Laplace equation.  
From the topological point of view, the Wiener test at $\infty$ presents thinness criteria
of sets near $\infty$ in fine topology. Precisely, the open set is a deleted neigborhood
of $\infty$ in fine topology if and only if $\infty$ is irregular.

{\bf Key words:} uniformly elliptic equations, metric compactification of $\rr^{N+1}$, 
measurable coefficients, regularity (or irregularity) of $\infty$, ${\cal A}$-harmonic measure, 
Dirichlet problem, PWB solution, ${\cal A}$-super- or subharmonicity, Newtonian capacity, Wiener test, fine topology, ${\cal A}$-thinness, Brownian motion, characteristic Markov process, differential generator

{\bf AMS subject classifications:} 35J25, 31C05, 31C15, 31C40, 60J45, 60J60

\newpage
\section{Description of Main Results}\label{description of results,historical remarks}
\subsection{Introduction and Motivation}\label{E:1:1}
{\large
This paper introduces the notion of regularity of  the point 
at infinity ($\infty$) and establishes a necessary and sufficient 
condition  for the unique solvability of the Dirichlet problem (DP)
in an arbitrary open set of $\rr^{N}$ for the uniformly elliptic 
equations in divergence form
\begin{equation}\label{Eq:W:1:1}
{\cal A}u=-(a_{ij}(x)u_{x_i})_{x_j}=0
\end{equation}
when the coefficients and boundary values 
are only supposed to be bounded and measuarable. 
The criterion is the same as Wiener test for the regularity
of $\infty$ concerning the classical DP for harmonic funtions 
in an arbitrary open set of $\rr^{N}$ \cite{Abdulla1}.

In order to formulate our result, we first introduce 
some terminology.
Let $\Omega \subset \rr^{N}$ $(N\ge3)$ denote any 
unbounded open 
subset, and $\pl \Omega$ its topological boundary.
We consider the differential operator ${\cal A}u$, with 
$a_{ij}=a_{ji}$ being real bounded measurable functions defined
in $\Omega$. 

Throughout this paper, we use the summation convention and assume that
${\cal A}$ is uniformly elliptic. That is, there is a constant 
$\lambda \ge 1$, such that
\begin{equation}\label{Eq:W:1:2}
\lambda^{-1}|\xi|^2 \le 
a_{ij}(x)\xi_i\xi_j\le \lambda |\xi|^2
\end{equation}
for all $x\in \Omega$ and all $\xi \in \rr^N$. We will assume that 
the coefficients $a_{ij}$ are defined and satisfy (\ref{Eq:W:1:2})
for all $x\in \rr^{N}$. This can always be achieved 
by putting $a_{ij}=\delta_{ij}$ outside of $\Omega$. 

Throughout, we use the standard notation of Sobolev spaces \cite{Adams}.
A function $u$ in $H_{loc}^{1,2}(\Omega)$ is a weak solution
of the equation (\ref{Eq:W:1:1}) in $\Omega$ if
\begin{equation}\label{Eq:W:1:3}
\int_\Omega a_{ij}u_{x_i}\phi_{x_j}dx=0
\end{equation}
whenever $\phi \in C_0^{\infty}(\Omega)$. A function $u$ in $H_{loc}^{1,2}(\Omega)$ is a supersolution
of (\ref{Eq:W:1:1}) in $\Omega$ if
\begin{equation}\label{Eq:W:1:4}
\int_\Omega a_{ij}u_{x_i}\phi_{x_j}dx\ge 0
\end{equation}
whenever $\phi \in C_0^{\infty}(\Omega)$ is nonnegative.
A function $u$ is a subsolution of (\ref{Eq:W:1:1}) in $\Omega$ if
$-u$ is a supersolution of (\ref{Eq:W:1:1}). 

The weak solution of (\ref{Eq:W:1:1})
is locally H\"{o}lder continuous \cite{De Giorgi, Nash, Moser}.
Continuous weak solution of (\ref{Eq:W:1:1}) in $\Omega$ 
is called ${\cal A}$-{\it harmonic} in $\Omega$. 

A function $u$ is called a ${\cal A}$-{\it superharmonic} in $\Omega$ if it satisfies the following conditions:
\begin{description}
\item{\bf(a)} $-\infty<u\le+\infty$, $u < +\infty$ on a dense subset of $\Omega$; 
\item{\bf(b)} $u$ is lower semicontinuous (l.s.c.);
\item{\bf(c)} for each open $U\Subset\Omega$ and each ${\cal A}$-harmonic $h\in C(\ov{\Omega})$,
the inequality $u\ge h$ on $\partial U$ implies $u \ge h$ in $U$.    
\end{description}
We use the notation
${\cal S}(\Omega)$ for a class of all ${\cal A}$-superharmonic functions in $\Omega$.
Similarly, $u$ is ${\cal A}$-{\it subharmonic} in $\Omega$ if $-u$ is 
${\cal A}$-superharmonic in $\Omega$; the class of all ${\cal A}$-subharmonic
functions in $\Omega$ is $-{\cal S}(\Omega)$. 

It is well known (\cite{LSW, HKM}) that $u\in H_{loc}^{1,2}$ is 
${\cal A}$-superharmonic if and only if
it is a supersolution with
\begin{equation}\label{Eq:W:1:5}
u(x)=ess \liminf_{y\to x} u(y) \ \quad\text{for all}\ x \in \Omega. 
\end{equation}
Moreover, in the $L_1$-equivalence
class of every supersolution, there is an ${\cal A}$-superharmonic representative
which satisfies (\ref{Eq:W:1:5}).

Given boundary function $f$ on $\partial \Omega$, consider a Dirichlet
problem (DP)
\begin{equation}\label{Eq:W:1:6}
{\cal A}u=0 \quad\text{in}\ \Omega, \ \ u=f \quad\text{on}\ \pl \Omega
\end{equation}

Consider a compactified space $\ov{\rr^{N}}=\rr^{N}\cup \{\infty\}$.
Denote by ${\bf \Omega}$ and ${\bf \pl\Omega}$ respectively, the metric 
compactification of $\Omega$ and $\pl\Omega$. Throughout, we always 
assume that $\infty \in {\bf \pl\Omega}$. 
For example, if $\Omega$ is an exterior of compact, then $\infty$ is an
isolated boundary point of ${\bf \Omega}$. 

Assuming for a moment
that $\pl \Omega$ is non-compact, $f\in C(\partial \Omega)$, and $f$ has a limit $f_o$ at infinity, prescribe $f(\infty)=f_o$. 
The {\it generalized upper (or lower) solution} (in the sense of Perron-Wiener-Brelot)  of the DP is defined as
\begin{equation}\label{Eq:W:1:7}
\ov H_f^{{\bf \Omega}} \equiv \inf \{u\in {\cal S}(\Omega): \liminf_{x\to y, x\in {\bf \Omega}} u \ge f(y) \ \quad\text{for all}\ y \in {\bf \pl \Omega} \}
\end{equation}
\begin{equation}\label{Eq:W:1:8}
\underline H_f^{{\bf \Omega}} \equiv \sup \{u\in {-\cal S}(\Omega): \limsup_{x\to y, x\in {\bf \Omega}} u \le f(y) \ \quad\text{for all}\ y \in {\bf \pl \Omega} \}
\end{equation}

According to classical theory (\cite{Doob, HKM}),  $f$ is a {\it resolutive boundary function}, in the sense that 
\[ \ov H_f^{{\bf \Omega}} \equiv \underline H_f^{{\bf \Omega}} \equiv  H_f^{{\bf \Omega}}. \]
Being ${\cal A}$-harmonic in $\Omega$, $H_f^{{\bf \Omega}}$ is called {\it a generalized solution} of the DP (\ref{Eq:W:1:5}).
The generalized solution is unique by construction, and it coincides with the classical, or
Hilbert, or Sobolev space solution whenever the latter exists. Also,
note that the construction of the generalized solution
is accomplished by prescribing the behavior of
the solution at $\infty$. The elegant theory, 
while identifying a class of unique solvability, 
leaves the following questions open:
\begin{itemize}
\item Would a unique solution still exist if its 
limit at infinity were not specified? That is, could it be 
that the solutions would pick up the ``boundary value" $f_o$
without being required?   
\item What if the boundary datum $f$ on $\pl \Omega$, while being continuous, 
does not have a limit at infinity, for example, it 
exhibits bounded oscillations. Is the DP uniquely solvable?
\end{itemize}
The answer to these fundamental questions depends on whether 
$\Omega$ is sufficiently sparse, or equivalently $\Omega^c$ is
sufficiently thin near $\infty$. In recent papers \cite{Abdulla1,Abdulla2,
Abdulla3} the answer is expressed in
terms of the Wiener test for the regularity of $\infty$
for Laplace and heat equations.
The principal purpose of this paper is to prove that 
the Wiener test for the regularity of $\infty$, and 
accordingly for the uniqueness of the bounded solution of the
DP, is the same as that for Laplace equation. This is the result 
in the spirit of \cite{LSW} concerning the regularity of finite 
boundary points.

\subsection{Formulation of Problems and Main Result}\label{E:1:2}

Furthermore, we assume that $f: \partial \Omega \rightarrow \rr$ 
is a {\it bounded Borel measurable} function. Bounded Borel measurable
functions are resolutive \cite{Doob}. 
Concerning $\Omega$, we don't exclude the case when $\pl \Omega$ is compact.
Without loss of generality, we only assume that
$\Omega$ has at least one connected component $\Omega_e$ such that 
$\infty \in {\bf \pl\Omega_e}$. Obviously if this assumption is not satisfied,
then the unbounded open set $\Omega$ is a union of (at most countable) connected bounded
components. In this case, within each bounded connected component, 
the solution is defined uniquely from boundary values on the boundary of this component. 
This makes the problem uniquely solvable in the whole region without prescribing
the boundary function at $\infty$.   

By fixing an arbitrary finite real number $\ov f$, extend a function
$f$ as $f(\infty)=\ov f$. Obviously, the extended function
is a bounded Borel measurable on ${\bf \pl\Omega}$. Since bounded
Borel measurable functions are resolutive (\cite{Doob}),
there exists a unique bounded generalized solution $H_f^{{\bf \Omega}}$.
It is natural to call it a generalized solution of the
DP (\ref{Eq:W:1:6}). The major question now becomes:\\

{\it {\bf Problem 1:} How many bounded solutions do we actually have, or does the constructed
solution depend on} $\ov f$ ? \\

Note that if in particular, 
$\pl \Omega$ is non-compact and $f$ has a limit
$f_o$ at $\infty$, the above constructed generalized solution is 
included here by choosing $\ov f = f_o$.

It is well possible that $H_f^{{\bf \Omega}}$ does not take continuously on the boundary values prescribed by $f$
at the finite boundary points. Therefore, the finite boundary point $w \in \pl \Omega$ is called {\it regular for} $\Omega$ if
\begin{equation}\label{Eq:W:1:9}
\lim_{z\to w, z\in \Omega} H_f^{{\bf \Omega}} (z) = f(w) \quad\text{for all bounded}\ f\in C(\pl \Omega).
\end{equation}
The regularity of a boundary point is a problem of local nature, and depends on the measure-geometric properties of the boundary in a neighborhood of the point, which is in turn dictated by the differential operator.  In his pioneering works \cite{Wiener1, Wiener2} Wiener discovered the
necessary and sufficient condition for the regularity of 
finite boundary points for harmonic functions. Remarkably, the regularity criteria
for finite boundary points with respect to the elliptic equation (\ref{Eq:W:1:1})
is the same \cite{LSW}.

For a given set $A$ from $\ov{\rr^{N}}$, denote as $1_A$ the indicator function
of $A$. Since the indicator function $1_\infty$ is a resolutive boundary function,
the ${\cal A}$-harmonic measure of $\infty$ is well defined (\cite{Doob, HKM}):

\[ \mu_{{\bf \Omega}}(\cdot,\infty)=H_{1_\infty}^{\bf \Omega}(\cdot). \]

It is said that $\infty$ is an ${\cal A}$-harmonic measure null set if $\mu_{{\bf \Omega}}(z,\infty)$
vanishes identically in ${\bf \Omega}$. If this is not the case, $\infty$ is a set of
positive ${\cal A}$-harmonic measure. Since $H_{1_\infty}^{\bf \Omega}$ is ${\cal A}$-harmonic
in $\Omega$, from the strong maximum principle it follows that if $\mu_{{\bf \Omega}}(x,\infty)>0$,
then it is positive in the whole connected component of $\Omega$ which contains $x$. 
We can now formulate the measure-theoretical counterpart
of the Problem 1:\\

{\it {\bf Problem 2:} Given $\Omega$, is the ${\cal A}$-harmonic measure 
of $\infty$ null or positive} ?\\

Note that the assumption $\infty\in{\bf \Omega_e}$ doesn't
cause a loss of generality in this context, and $\infty$ is an ${\cal A}$-harmonic measure 
null set otherwise.

In fact, both major problems are equivalent, and the next definition expresses the 
connection between them.\\  

{\bf Definition 1.1.}\ $\infty$ is said to be {\it regular} for
${\bf \Omega}$ if it is an ${\cal A}$-harmonic measure null set. Conversely, $\infty$ {\it is irregular} if it has a positive ${\cal A}$-harmonic measure.\\

The notion of the regularity of $\infty$ is, in particular, related to 
the notion of continuity of the solution at $\infty$. 

{\it {\bf Problem 3:} Given $\Omega$ with non-compact $\pl\Omega$,
whether or not}

\[\liminf f \le \liminf (H_f^{{\bf \Omega}}) \le \limsup (H_f^{{\bf \Omega}}) \le \limsup f \]
\begin{equation}\label{Eq:W:1:10}
 \quad\text{as}\ z \rightarrow \infty \quad\text{for all bounded}\ f\in C(\pl \Omega).
\end{equation}

Note that if $f$ has a limit at $\infty$, (\ref{Eq:W:1:10}) simply means that 
the solution $H_f^{{\bf \Omega}}$ is continuous at $\infty$.\\

The notion of the regularity of $\infty$ introduced in Definition 1.1 and earlier in \cite{Abdulla1, Abdulla2, Abdulla3} fits naturally in the framework of ${\cal A}$-fine topology. ${\cal A}$-fine topology
is the coarsest topology of $\rr^N$ which makes every ${\cal A}$-superharmonic function continuous
\cite{Doob, HKM}.  ${\cal A}$-fine topology is finer than the Euclidean topology.
A major problem is to find the structure of the neigborhood base in
${\cal A}$-fine topology. It is well-known that there is an elegant connection
between this problem and the problem on the regularity of finite boundary points.
Namely, given open set $\Omega \subset \rr^N$, its finite boundary point $x_0$
is irregular if and only if $\Omega$ is a deleted neigborhood of $x_0$
in fine topology \cite{Doob, Landkov, HKM}. Our definition of the 
regularity of $\infty$ reveals a similar connection
for the point at infinity.\\

Let $u: \rr^N \rightarrow (-\infty, +\infty]$ be an arbitrary ${\cal A}$-superharmonic function.
Extend it to $\ov{\rr^{N}}$ as follows:
\[ u(\infty)=\liminf_{x\to \infty} u(x). \]
What is the coarsest topology in $\ov{\rr^{N}}$ which makes every extended ${\cal A}$-superharmonic function continuous? How should ${\cal A}$-fine topology of $\rr^N$ be extended
to $\infty$ with ``minimum increase" of the Euclidean neigborhood base of $\infty$, to guarantee
continuity at $\infty$ of every extended ${\cal A}$-superharmonic function? 
The following definition is helpful in understanding the structure of the neigborhood
base of $\infty$ in ${\cal A}$-fine topology.\\

{\bf Definition 1.2.}\ A subset $E$ of $\rr^N$ is called ${\cal A}$-thin at $\infty$
in the following two cases:
\begin{description}
\item{\bf(a)} $E$ is bounded
\item{\bf(b)} $E$ is unbounded and there exists an ${\cal A}$-superharmonic 
function $u$ in $\rr^N$ such that
\begin{equation}\label{Eq:W:1:11}
 u(\infty)<\liminf_{x\to\infty, x\in E}u(x).
\end{equation}
\end{description}

A set $E$ is ${\cal A}$-thin at $\infty$ if and only if $E^c$ is a deleted ${\cal A}$-fine
neigborhood of $\infty$ (see Lemma ~\ref{4}). Furthermore, the ${\cal A}$-fine derived set
of $E$, that is, the ${\cal A}$-fine closed set of ${\cal A}$-fine limit points of $E$,
will be denoted by $E^f$. Hence, $E$ is ${\cal A}$-thin at $\infty$ if and only if
$\infty \not\in E^f$.  As in the case of the finite points, ${\cal A}$-fine
topology
is strictly finer than Euclidean topology near $\infty$. The exterior of any compact set is a deleted neigborhood of $\infty$ both
in ${\cal A}$-fine and Euclidean topology. However, in ${\cal A}$-fine topology, there are deleted neigborhoods
of $\infty$, which are unbounded open sets with non-compact boundaries. 

We can now formulate the topological counterpart of problems 1, 2 and 3:\\

{\it {\bf Problem 4:} Is given open set ${\bf \Omega}$ a deleted neigborhood of $\infty$
in ${\cal A}$-fine topology? Conversely, is $\Omega^c$ ${\cal A}$-thin at $\infty$?} \\

The principal result of this paper expresses the solutions to equivalent
Problems 1-4 in terms of the Wiener test for the regularity of $\infty$. 

Recall that if $K\subset \rr^N$ is compact, the Newtonian capacity of $K$ is
\[ cap(K)\equiv sup\{\mu(\rr^N): \mu \in {\cal M}(K), V_\mu \le 1 \quad\text{in}\ \rr^N \}, \]
where ${\cal M}(K)$ denotes the collection of all nonnegative Radon measures on $\rr^N$ with
support in $K$, and
\[ V_\mu(x) \equiv \int_{\rr^N} \frac{d\mu(y)}{|x-y|^{N-2}}, \ x\in \rr^N  \]
is the Newtonian potential of $\mu$. Let
\[ \Gamma_n\equiv cap(E_n), \ E_n=\Omega^c\cap\{x:2^{n-1}\le|x|\le2^n \}  \]
Our main theorem reads:
\begin{theorem}\label{main theorem}
The following conditions are equivalent:
\begin{description}
\item{(1)\ } $\infty$ is regular (or irregular)
\item{(2)\ }DP has a unique (or infinitely many) bounded solution(s)
\item{(3)\ } If $\pl \Omega$ is non-compact, (\ref{Eq:W:1:10}) is satisfied (respectively isn't satisfied).
\item{(4)\ } Wiener series
\begin{equation}\label{Eq:W:1:12}
\sum_{n} 2^{-n(N-2)}\Gamma_n
\end{equation}
diverges or converges.
\item{(5)\ } $\infty \in(\Omega^c)^f$ (or $\Omega^c$ is ${\cal A}$-thin at $\infty$)\end{description}
\end{theorem}

An equivalent characterization is valid if one replaces (\ref{Eq:W:1:12})
with the series
\[\sum_{n} \lambda^{-n}\gamma_n \]
where $\lambda >1$, $\gamma_n=cap(e_n), e_n=\Omega^c \cap \{x:\lambda^{-n}\le\Phi(x)\le\lambda^{-n+1}\}$ and 
\[ \Phi(x)=\frac{\Gamma(\frac{N}{2})}{(N-2)(4\pi)^{n/2}}|x|^{2-N}  \]
is the fundamental solution of the Laplace equation. Another equivalent characterization is valid if one replaces the series (\ref{Eq:W:1:12}) with the Wiener integral
\[ \int_1^{+\infty} \frac{c(\rho)}{\rho^{N-1}}d\rho  \]
where
\[ c(\rho)\equiv cap(\Omega^c\cap\{x:|x|\le \rho\}). \]
\\
{\bf Counterpart at the finite boundary point:} In a paper \cite{SerrinW} a generalization
of the well-known Kelvin transformation for Laplace equation was introduced
for the equation (\ref{Eq:W:1:1}). This transformation allows to transform the uniformly elliptic equation (\ref{Eq:W:1:1}) near $\infty$ to another uniformly elliptic equation
near finite boundary point $x_0$ which is the inversion of $\infty$. Moreover, ${\cal A}$-harmonic functions which are bounded near $\infty$ will be transformed to $\ov{{\cal A}}$-harmonic functions of class
\begin{equation}\label{Eq:W:1:13}
u(x)=O(|x-x_0|^{2-N}), \ \quad\text{as}\ x\rightarrow x_0
\end{equation}
where $\ov{{\cal A}}$ is an elliptic operator of type (\ref{Eq:W:1:1}),(\ref{Eq:W:1:2}).
Counterpart of the {\bf Problem 1} is whether or not  singular Dirichlet problem in a bounded open set
$\Omega$, with $x_0\in \partial\Omega$, is uniquely solvable in a class (\ref{Eq:W:1:13}).
Otherwise speaking, whether or not fundamental solution kind singularity is removable
in a non-isolated boundary point $x_0$. Since divergence or convergence of the Wiener series is preserved under the inversion transformation, Theorem~\ref{main theorem} implies that
the classical Wiener test at $x_0$ is a necessary and sufficient condition 
for uniqueness in a singular Dirichlet problem, as well as for the removability 
of the singularity according to (\ref{Eq:W:1:13}).   
\\
{\bf Probabilistic Counterpart:} From the probabilistic standpoint, Wiener test 
of Theorem ~\ref{main theorem} presents an asymptotic probability law for Markov processes
in $\rr^N$. Assume that $\{X_t\}$ is a continuous time, time-homogeneous Markov
process with infinitesimal Dynkin generator being a differential operator $-{\cal A}$\cite{Dynkin}.
Conversely, $X_t$ is the characteristic process for the differential operator   
$-{\cal A}$. The characteristic process of the operator $-{\cal A}\equiv \frac{1}{2}\Delta$
is the Wiener process or Brownian motion. Assume that, $B$ is a closed set in $\rr^N$ clustering to
$\infty$, and $B_n$ is the intersection of $B$ with the spherical shell $2^{n-1}\le |x| < 2^n$.
Let ${\bf B}$ be an event that $\{t: X_t\in B\}$ clusters to $+\infty$. Then
\[ P_{\bullet}({\bf B})=0 \quad\text{or}\ 1 \quad\text{according as}\ \sum_n 2^{-n(N-2)} cap(B_n) < \quad\text{or}\ = +\infty.  \]
The same result for an $N$-dimensional standard random walk on the $N$-dimensional lattice of points with
integer coordinates was proved in \cite{ito1}. In \cite{ito2} (p.257) the law is mentioned 
for standard $N$-dimensional Brownian motion. 
\subsection{Historical Comments}\label{E:1:2}
It will be convenient to make some remarks concerning the Dirichlet problem for uniformly elliptic equations.  The solvability, in some generalized sense, of the classical DP in bounded open set $E \subset \rr^{N}$, with prescribed data on $\partial E$, is realized within the class of resolutive boundary functions, identified by Perron's method and its Wiener \cite{Wiener1,Wiener2} and Brelot \cite{Brelot} refinements.
Such a method is referred to as the PWB method, and the corresponding solutions are PWB solutions. The main tool of the PWB method consists of Harnack estimates for the 
Laplace equation. 

Wiener, in his pioneering works \cite{Wiener1,Wiener2}, proved a necessary and sufficient condition for the finite boundary point $x_o \in \partial E$ to be regular in terms of the ``thinness" of the complementary set in the neighborhood of $x_o$. Cartan pointed out that the thinness could be characterized by means of
fine topology -- the coarsest topology of $\rr^N$ which makes every superharmonic function continuous.
In fact, a finite boundary point is irregular if and only if $E$ is a deleted neigborhood 
of $x_o$ in fine topology \cite{Doob}. 

De Giorgi \cite{De Giorgi} and Nash \cite{Nash} almost simultaneously 
proved that any local weak solution of (\ref{Eq:W:1:1}) is locally H\"older 
continuous. Moser \cite{Moser} gave a simpler proof of this fact as well as the
Harnack inequality. In \cite{LSW} it is proved that the Wiener test for the regularity
of finite boundary points with respect to elliptic operator (\ref{Eq:W:1:1})
coincides with the classical Wiener test for the boundary regularity of harmonic functions.
Hence, the fine-topological neigborhood base of the finite boundary point is 
independent of elliptic operator (\ref{Eq:W:1:1}). The Wiener test for the regularity of finite boundary points for linear degenerate elliptic equations is proved in \cite{FJK}. 
Wiener test for the regularity of finite boundary points for quasilinear 
elliptic equations was settled due to \cite{Mazya3, gariepy, martio, kilpi}.
Nonlinear potential theory was developed along the same lines as classical potential theory
for the Laplace operator and we refer to monographs \cite{LSW, Ziemer}.

Despite the importance of the Wiener criterion in Analysis, its meaning 
for the point at infinity was not correctly understood. In the framework of  
Brelot's theory and its generalizations, the regularity of $\infty$ was associated with the existence of the solution to the Dirichlet problem with the same limit at $\infty$ as the boundary function. 
This approach implied that $\infty$ is always regular if $N\ge 3$, and acordingly, the geometric nature of the Wiener criterion was ignored. 
``Labeling" of $\infty$ as ``always regular" became a standard
result in classical and modern potential theory.
The deficiency of this approach is clear in the context of fine topology.
The connection mentioned above between the irregularity of the boundary point and 
fine topological thinness falls apart for the point at infinity, since
otherwise fine topology would be trivial at $\infty$. Another
deficiency comes out in the context of the asymptotic properties
of Markov processes. A well known elegant connection
between the irregularity of boundary points for the 
Laplacian and non-escaping property of Wiener processes
starting at boundary point is ignored for the point $\infty$.  

In \cite{Abdulla1,Abdulla2,Abdulla3} we introduced a correct notion of regularity 
of the point at infinity on $\pl E$ for the Laplace and heat equations.
Basically, the 
Dirichlet problem for $E$ with continuous data 
$\phi$, has either one and only one bounded solution, or 
infinitely many. If the DP has a unique solution,  
the point at infinity on $\pl E$ is regular, otherwise 
it is irregular. In \cite{Abdulla1} we characterize the regularity of $\infty$
through the Wiener test at $\infty$.   
The principal result of this paper proves that the Wiener test
for the regularity of $\infty$ is independent 
of elliptic operators (\ref{Eq:W:1:1}) with bounded measurable coefficients.
The Wiener test at $\infty$ characterizes fine topological thinness
at $\infty$. It also characterizes the asymptotic properties of the
characteristic Markov processes with differential generator $-{\cal A}$.

\section{Preliminary Results on Potential Theory}\label{preliminaries}

This section formulates basic known facts about ${\cal A}$-capacity,
${\cal A}$ potentials and Green functions for differential operator (\ref{Eq:W:1:1})
which we need to prove the main Theorem ~\ref{main theorem}. 
Lemmas ~\ref{1} -- ~\ref{3} are due to \cite{LSW}. In Lemma ~\ref{4}, we formulate
${\cal A}$-thinness criteria at $\infty$ in ${\cal A}$-fine topology. 

Let $\Sigma$ be a fixed open sphere and $E$ be a compact subset of $\Sigma$.
Then the ${\cal A}$-capacity of $E$ with respect to operator ${\cal A}$ and sphere $\Sigma$
is defined as
\[ cap_{\cal A}(E)\equiv \inf \Big \{ D_{\cal A} : \phi\in H_0^{1,2}(\Sigma), \phi \ge 1 \quad\text{on}\ E \quad\text{in the sense of}\ H_0^{1,2}(\Sigma) \Big \}  \]
where
\[ D_{\cal A} = \int_\Sigma a_{ij} \phi_{x_i} \phi_{x_j} dx \]
The function $u$ giving the infimum to $D_{\cal A}(\phi)$ is called ${\cal A}$-capacitary
potential (with respect to ${\cal A}$ and $\Sigma$). 

\begin{lemma}\label{1}
\begin{description}
\item{(1)\ } There exists one and only one ${\cal A}$-capacitary potential $u\in H_0^{1,2}(\Sigma)$ such that
$u= 1$ on $E$ in the sense of $H_0^{1,2}(\Sigma)$ and $cap_{\cal A}(E)=D_{\cal A}(u)$.
\item{(2)\ } ${\cal A}$-capacitary potential is ${\cal A}$-supersolution in $\Sigma$, and ${\cal A}$-harmonic
in $\Sigma-E$.
\end{description}
\end{lemma}

For any Radon measure $\mu$ with compact support in $\Sigma$, consider a problem
\begin{equation}\label{Eq:W:2:1}
{\cal A}u=\mu \quad\text{in}\ \Sigma,\ \  u=0 \quad\text{on}\ \pl\Sigma.
\end{equation}

$u\in L^1(\Sigma)$ is a weak solution of (\ref{Eq:W:2:1}) if
\[ \int_\Sigma u {\cal A}\Phi dx = \int_\Sigma \Phi d\mu  \]
for every $\Phi\in H_0^{1,2}(\Sigma)\cap C(\ov{\Sigma})$ such that ${\cal A}\Phi \in C(\ov{\Sigma})$.

\begin{lemma}\label{2}
\begin{description}
\item{(1)\ } There exists a unique weak solution $u$ of (\ref{Eq:W:2:1}) such that $u\in H_0^{1,p}(\Sigma)$ for every $p<n/(n-1)$.
\item{(2)\ }  Weak solution $u\in H_0^{1,2}(\Sigma)$, if and only if $\mu \in H^{-1,2}(\Sigma)$. 
\end{description}
\end{lemma} 

The Green's function $g(x,y)$ of the operator ${\cal A}$ on $\Sigma$ is defined as the weak solution 
of the problem (\ref{Eq:W:2:1}) with $\mu=\delta_y$, where $\delta_y$ is the Dirac measure of $y$.

\begin{lemma}\label{3}
\begin{description}
\item{(1)\ } For every Radon measure with compact support in $\Sigma$, the integral
\begin{equation}\label{Eq:W:2:2}
u(x)=\int_\Sigma g(x,y)d\mu(y).
\end{equation}
exists and finite a.e., and is a weak solution of (\ref{Eq:W:2:1}). 
\item{(2)\ }  ${\cal A}$-capacitary potential is a weak solution of (\ref{Eq:W:2:1}) with the
nonnegative (${\cal A}$-capacitary) measure supported on the exterior boundary of compact $E$,
and accordingly the representation (\ref{Eq:W:2:2}) is valid. Moreover, $\mu(E)=cap_{\cal A}(E)$.
\item{(3)\ } Let $g$ and $\ov{g}$  be the Green functions for any uniformly elliptic operators ${\cal A}$ and $\ov{{\cal A}}$ with the  ellipticity constant $\lambda$ on a sphere $\Sigma$. Then, for any compact subset $E$ of $\Sigma$, there exists a constant $K$ depending on $E, \Sigma$ and $\lambda$ such that
\begin{equation}\label{Eq:W:2:3}
K^{-1}\ov{g}(x,y) \le g(x,y) \le K \ov{g}(x,y), \quad\text{for all}\ x,y\in \Sigma.
\end{equation}
\item{(4)\ } If $\mu$ is a nonnegative  Radon measure with compact support $E$ in $\Sigma$,
and $u$ and $\ov{u}$ are the weak solutions of (\ref{Eq:W:2:1}) for differential operators 
${\cal A}$ and $\ov{{\cal A}}$ respectively, then
\begin{equation}\label{Eq:W:2:4}
\lambda^{-2}cap_{\ov{\cal A}}(E) \le cap_{\cal A}(E) \le \lambda^2cap_{\ov{\cal A}}(E) 
\end{equation}
\begin{equation}\label{Eq:W:2:5}
K^{-1}\ov{u}(x) \le u(x) \le K \ov{u}(x), \quad\text{a.e. on}\ E
\end{equation}
If $\mu \in H^{-1,2}(\Sigma)$, then (\ref{Eq:W:2:5}) is valid in the sense of $H_0^{1,2}(\Sigma)$.
\item{(5)\ } $cap_{\cal A}(\cdot)$ is a monotone, subadditive set function; it is a Choquet 
capacity (\cite{Doob, HKM}) and for any differential operators  ${\cal A}$ and $\ov{{\cal A}}$, $cap_{\cal A}(\cdot)$ and 
$cap_{\ov{{\cal A}}}(\cdot)$ are mutually absolutely continuous (see (\ref{Eq:W:2:4})).
\item{(6)\ } In the $L_1$-equivalence
class of every capacitary potential, there is an ${\cal A}$-superharmonic representative
which satisfies (\ref{Eq:W:1:5}).
\item{(7)\ } $g(\cdot,y)\ge 0$ is ${\cal A}$-harmonic and H\"{o}lder continuous
in $\Sigma-y$, $\lim_{x\to y} g(x,y)=+\infty$. 
\end{description}
\end{lemma} 

If $\ov{{\cal A}}$ is taken to be the Laplace operator, then from (\ref{Eq:W:2:3})
it follows that  for any compact subset $E$ of the open sphere $\Sigma$, 
\begin{equation}\label{Eq:W:2:6}
K^{-1}|x-y|^{2-N} \le g(x,y) \le K |x-y|^{2-N} \quad\text{for all}\ x,y\in E
\end{equation}
where $K$ depends on $E, \Sigma, \lambda$ and $R$ is a radius of $\Sigma$.

If $N\ge3$ and the radius of $\Sigma$ goes to infinity, the Green function increases
to a Green function $G(x,y)$ of $\rr^N$, which is locally H\"{o}lder continuous and 
${\cal A}$-harmonic in $\rr^N-y$. In fact, $G(\cdot,y)\in H_{loc}^{1,p}(\rr^N) \cap H_{loc}^{1,2}(\rr^N-y), p<n/(n-1)$. Moreover, we have
\begin{equation}\label{Eq:W:2:8}
K^{-1}|x-y|^{2-N} \le G(x,y) \le K |x-y|^{2-N}, 
\end{equation}
uniformly in $\rr^N$. Accordingly, ${\cal A}$-capacity and related ${\cal A}$-capacitary potential 
are well defined if $\Sigma$ and $g(x,t)$ are replaced with $\rr^N$ and $G(x,t)$ respectively. 

In particular, the ${\cal A}$-capacitary 
potential $\zeta_E$ of the compact subset $E\subset\rr^N$ is 
a weak solution of the problem (\ref{Eq:W:2:1}) with $\Sigma$ replaced 
by $\rr^N$. Accordingly, we have
\begin{equation}\label{Eq:W:2:9}
\zeta_E(x)=\int_{\rr^N} G(x,y)d\mu(y).
\end{equation} 
where $\mu$ is a ${\cal A}$-capacitary measure with support on the exterior boundary of 
$E$ with $\mu(\rr^N)=cap_{\cal A}(E)$. As before, in the $L_1$-equivalence class,
there is a lower semicontinuous representative which satisfies
(\ref{Eq:W:1:5}), and it is ${\cal A}$-superharmonic in $\rr^N$. 
Using different terminology, this representative is a smoothed ${\cal A}$-reduction
of $1$ on $E$. 

Furthermore, we will assume that the ${\cal A}$-capacity and ${\cal A}$-potential
of the compact subsets of $\rr^N$ are defined with respect to $\rr^N$.
We will drop subscript ${\cal A}$ when differential operator ${\cal A}$ is Laplacian.
In this case ${\cal A}$-capacity and ${\cal A}$-potential coincide with 
Newtonian capacity and Newtonian potential. Hence, in view of (\ref{Eq:W:2:8}) and (\ref{Eq:W:2:9}), (\ref{Eq:W:2:5}) becomes  
\begin{equation}\label{Eq:W:2:10}
K^{-1} V_{\nu}(x) \le \zeta_E(x) \le K V_{\nu}(x), 
\end{equation}
a.e. on $E$ and everywhere on $E^c$ with $\nu$ being the Newtonian
capacitary measure of $E$. 

Throughout, we will say that a property holds {\it quasieverywhere}, if it holds 
except on a set of ${\cal A}$-capacity (or Newtonian capacity) zero. Sets 
of capacity zero are called ${\cal A}$-{\it polar sets} (or simply
{\it polar sets}) in potential theory. Polar sets are essential for the description of the singularities
of ${\cal A}$-superharmonic functions. In fact, for any polar set $E \in \rr^N$,
there is an ${\cal A}$-superharmonic function $u$ in $\rr^N$ such that $u=+\infty$
on $E$. We refer to \cite{Doob, HKM} for the essential properties of polar sets.

Let $D_e$ be the connected component of $E^c$ which contains
$\infty$ and $E_e=\pl D_e$. The restriction of $\zeta_E$ to $D_e$ is ${\cal A}$-harmonic
and solves the exterior Dirichlet problem with boundary value $1$ on $E_e$ and
$0$ at infinity. More precisely, $\zeta_E \equiv H_{1_E}^{{\bf E^c}}$. From (\ref{Eq:W:2:8}),(\ref{Eq:W:2:9}) it follows that it
vanishes at infinity with the same rate as a fundamental solution of the Laplacian.
For example, if $E=\Sigma\equiv \{x:|x|<R\}$ then
\begin{equation}\label{Eq:W:2:11}
K^{-1} R^{N-2} |x|^{N-2} \le \zeta_\Sigma(x)\le K R^{N-2} |x|^{2-N} \quad\text{for}\  |x|>R
\end{equation}
Remembering that the Newtonian capacity of $\Sigma$ is $R^{N-2}$, we also have   
\begin{equation}\label{Eq:W:2:12}
\lambda^{-2}R^{N-2}\le cap_{\cal A}(\Sigma)\le \lambda^2 R^{N-2}
\end{equation}  

In the next lemma, we formulate the ${\cal A}$-thinness criteria of sets at $\infty$.

\begin{lemma}\label{4}
\begin{description}
\item{(1)\ } If the sets $E_1,E_2,...,E_n$ are ${\cal A}$-thin at $\infty$, then $E=\cup_{i=1}^n E_i$ is also ${\cal A}$-thin at $\infty$.
\item{(2)\ } A set $E$ is ${\cal A}$-thin at $\infty$ if and only if $E^c$ is a deleted
neigborhood of $\infty$ in ${\cal A}$-fine topology.
\end{description}
\end{lemma}

The proof of these facts is standard (\cite{Doob, Landkov, HKM}). For (1), it is easy to observe that if $v_i(x)$ is the ${\cal A}$-superharmonic function related to $E_i$ by the relation
(\ref{Eq:W:1:11}), then
\[ \sum_{i=1}^n v_i(x)  \]
will, by the lower semicontinuity of $v_i$, satisfy (\ref{Eq:W:1:11}).

To prove (2), first note that ${\cal A}$-fine topology is the coarsest 
topology where the sets of the form $\{u>\beta\}$, $\{u<\beta\}$ are open,
for all ${\cal A}$-superharmonic functions $u$ and for all real numbers $\beta$.
The family formed by the finite intersections of these sets forms a neigborhood base
of  ${\cal A}$-fine topology. Since $u$ is lower semicontinuous, the 
sets of the form $\{u>\beta\}$ are open in Euclidean topology. 
Accordingly, we can consider a neigborhood base of $\infty$ consisting 
of the sets 
\begin{equation}\label{Eq:W:2:13}
\cap_{i=1}^m \{x\in \Sigma^c: u_i \le \beta \}
\end{equation}
where $m$ is an integer, $\beta>0$, $\Sigma$ is an open sphere in $\rr^N$, $u_i$
is a locally bounded ${\cal A}$-superharmonic function in $\rr^N$ such that
$u_i(\infty)=0$. 

Assume that $E^c$ is a deleted ${\cal A}$-fine neigborhood of $\infty$.
Then it contains some element of the neigborhood base:
\[ \cap_{i=1}^m \{x\in \Sigma^c: u_i \le \beta \} \subset E^c \]
which means that
\[ E \subset \cup_{i=1}^m E_i \cup \Sigma \]
where $E_i=\{u_i>\beta \}$. By  Definition 2.1, both $\Sigma$ and each of the
sets $E_i$ are ${\cal A}$-thin at $\infty$. By assertion (1) of Lemma~\ref{4},
the union, and accordingly its subset $E$, is ${\cal A}$-thin at $\infty$.

Assume that $E$ is ${\cal A}$-thin at $\infty$. If $E$ is bounded in $\rr^N$, then
the assertion of the lemma is trivial. Assume $E$ is unbounded and choose an
${\cal A}$-superharmonic function in $\rr^N$ according to Definition 2.1
and a real number $\beta$ such that
\[ u(\infty)<\beta<\liminf_{x\to\infty, x\in E}u(x)  \]
Accordingly, there is an open sphere $\Sigma$ such that
\[ u(x)\ge \beta \quad\text{for all}\ x \in E\cap \Sigma^c  \]
This implies that ${\cal A}$-fine neigborhood $\{x\in \Sigma^c: u(x)<\beta \}$
of $\infty$ is contained in $E^c$. Hence, $E^c$ is a deleted ${\cal A}$-fine
neigborhood of $\infty$. 
 
\section{Proof of Theorem 1.1}\label{proof of the main theorem}
${\bf (1)\Leftrightarrow(2)}$: Assume that $\infty$ is regular and let $u_1$ and $u_2$ be 
two bounded solutions of DP. Then $v=u_1-u_2$
is a bounded solution of DP with zero boundary data on finite boundary points. 
Since a generalized solution is order preserving (\cite{Doob, HKM}), we have

\[ |v|\le H_{M\cdot1_\infty}^{\bf \Omega} \equiv M H_{1_\infty}^{\bf \Omega} \equiv 0, \quad\text{with}\> M=\sup|v|. \]
On the other hand, if $\infty$ is irregular, then for an arbitrary real number $r$, $r H_{1_\infty}^{\bf \Omega}$ is a generalized solution of the DP with zero boundary data on $\pl \Omega$.\\

${\bf (2)\Rightarrow(3)}$: Assume that (3) is not satisfied. That is to say, there is a bounded function $f$ such that
for some generalized solution $H_f^{\bf \Omega}$, 
one of the inequalities in (\ref{Eq:W:1:10}) is violated. In this case, by choosing a number $\ov{f}$ satisfying 

\[ f_{*}\equiv \liminf_{z\to\infty} f(z) \le \ov{f} \le \limsup_{z\to\infty} f(z)\equiv f^{*}, \]
and by extending $f(\infty)=\ov{f}$, we can always construct a generalized solution $H_f^{\bf \Omega}$ which satisfies (\ref{Eq:W:1:10}) . Indeed, since a generalized solution is order preserving, we clearly have
\begin{equation}\label{Eq:W:3:1}
|H_f^{\bf \Omega}| \le M \equiv \sup |f|.
\end{equation}
Then for an arbitrary $\epsilon > 0$ we choose $R>0$ such that 
\begin{equation}\label{Eq:W:3:2}
f_* - \epsilon \leq f \leq f^* + \epsilon, \quad\text{ on }\> \pl \Omega \cap \Sigma^c,
\end{equation}
where $\Sigma \equiv \{ z \mid |z| < R \}$.
Since capacitary potential $\zeta_\Sigma \equiv H_{1_{\pl\Sigma}}^{{\bf \Sigma^c}}$ we have 
\begin{equation}\label{Eq:W:3:3}
f_* - \epsilon - 2M\zeta_\Sigma \leq H_f^{\bf \Omega} \leq  f^* + \epsilon + 2M\zeta_\Sigma \quad\text{ on }\> {\bf \Omega \cap \Sigma^c}.
\end{equation}
This follows from the fact that (\ref{Eq:W:3:3}) is satisfied on ${\bf \pl (\Omega\cap \Sigma^c)}$ and a generalized
solution is order preserving. Passing to limit, first as $x\rightarrow \infty$, and then as $\epsilon \downarrow 0$, from (\ref{Eq:W:3:3}) and (\ref{Eq:W:2:11}) it follows that the constructed generalized solution satisfies  (\ref{Eq:W:1:10}). Contradiction with uniqueness.\\  

${\bf (3)\Rightarrow(2)}$:  Let $u_1$ and $u_2$ be two bounded solutions of the DP. Their
difference is a generalized solution with a zero boundary function, and accordingly,
it vanishes identically in view of (\ref{Eq:W:1:10}). On the other hand, if at least for one $f$ (\ref{Eq:W:1:10})
is violated, then from the relation $(2)\Rightarrow(3)$ for direct assertions, it follows that there are
at least two solutions of the DP. This implies that the DP with zero boundary data on $\pl \Omega$ has a
non-trivial solution. That means it has infinitely many solutions. \\ 

${\bf (1)\Rightarrow(4)}$:  Assume that the series (\ref{Eq:W:1:12}) converges, and prove that $H_{1_\infty}^{\bf \Omega} \not\equiv 0$. In fact, we
are going to prove that
\begin{equation}\label{Eq:W:3:4}
\limsup_{z\to\infty}H_{1_\infty}^{\bf \Omega}=1.
\end{equation}
Let $0<\ep<1$ be an arbitrary small number. Choose a positive integer $m$ so large that
\begin{equation}\label{Eq:W:3:5}
\sum_{n=m+1}^{\infty} 2^{-n(N-2)}\Gamma_n < \frac{2\ep}{K4^{N-1}}.
\end{equation}
Since $\Omega$ has at least one connected component $\Omega_e$ such that 
$\infty \in {\bf \pl\Omega_e}$, we can choose $m$ so large that there is a point $x^*$ with
\[ x^* \in \Omega_e \cap \ov{\Sigma}_{2^{m-1}}  \]
Consider a sequence of increasing ${\cal A}$-harmonic functions
\[ \vartheta_M(x)=\sum_{n=m+1}^{M} \zeta_{E_n}(x), M=m+1,m+2,... \ x\in \Omega^m \equiv \Big ( 
\bigcup_{n=m+1}^{+\infty} E_n \Big )^c, \]
where $\zeta_{E_n}$ is a capacitary potential of $E_n$. 
Obviously we have
\[ \vartheta(x) \equiv \sum_{n=m+1}^{+\infty}\zeta_{E_n}(x)\equiv \sup_{M}\vartheta_M(x), \ x\in \Omega_m \]
Since $\vartheta$ is a limit of the increasing sequence of ${\cal A}$-harmonic
functions in each component of $\Omega^m$ which includes $\Omega$, either
$\vartheta\equiv +\infty$, or $\vartheta$ is ${\cal A}$-harmonic. Let us 
now estimate $\vartheta(x^*)$. From (\ref{Eq:W:2:8})-(\ref{Eq:W:2:10})
it follows that
\[ \zeta_{E_n}(x^*)\le K \int_{E_n} |y-x^*|^{2-N}d\mu_n(y) \le K(2^{n-1}-|x^*|)^{2-N} \Gamma_n \]
\[ \le K 4^{N-2} 2^{-n(N-2)}\Gamma_n, 
\quad\text{ for }\ n\ge m+1, \]
\[ \vartheta(x^*)\le K 4^{N-2} \sum_{n=m+1}^{+\infty} 2^{-n(N-2)}\Gamma_n \]
and in view of (\ref{Eq:W:3:5}) we have
\begin{equation}\label{Eq:W:3:6}
\vartheta(x^*)<\frac{\ep}{2}
\end{equation}
Hence, $\vartheta$  is ${\cal A}$-harmonic in $\Omega_m$ and ${\cal A}$-superharmonic
in $\rr^N$.  In fact, $\vartheta$ is a generalized
solution of the DP in $\Omega^m$ under the boundary function
\begin{equation}\label{Eq:W:3:7}
\sum_{n=m+1}^{+\infty} \zeta_{E_n}(x)\ge 1 \quad\text{ quasieverywhere on }\ \ \pl\Omega^m.
\end{equation} 
Since the generalized solution is order preserving, we have 
\begin{equation}\label{Eq:W:3:8}
0\le H_{1-1_\infty}^{\bf \Omega}\le \vartheta + \zeta_{\Sigma_{2^m}}  \quad\text{ on }\> \ \Omega\cap\Sigma_{2^m}^c.
\end{equation}
Let $R$ be an arbitrary number satisfying 
\begin{equation}\label{Eq:W:3:9}
R>2^m K^{\frac{1}{N-2}} \Big ( \frac{\ep}{2} \Big )^{\frac{1}{2-N}}
\end{equation}
Since $\vartheta$ is ${\cal A}$-harmonic in $\Omega^m$, from (\ref{Eq:W:3:6}) and
(\ref{Eq:W:3:7}) we conclude 
that there must be a point $x_R\in \Omega_e\cap\{x:|x|=R\}$ such that
\[ \vartheta(x_R)<\frac{\ep}{2} \] 
From (\ref{Eq:W:3:8}),(\ref{Eq:W:3:9}) and (\ref{Eq:W:2:11}) it follows that  
\begin{equation}\label{Eq:W:3:10}
0\le H_{1-1_\infty}^{\bf \Omega}(x_R) \le \ep.
\end{equation}
Passing to the limit, first as $R\rightarrow +\infty$, and then as $\ep \rightarrow 0$ from  (\ref{Eq:W:3:10}) it follows that 
\begin{equation}\label{Eq:W:3:11}
\liminf_{x\to\infty} H_{1-1_\infty}^{\bf \Omega} =0.
\end{equation}
Since
\[ H_{1-1_\infty}^{\bf \Omega}\equiv 1-H_{1_\infty}^{\bf \Omega}, \]
we arrive at (\ref{Eq:W:3:4}).\\

${\bf (4)\Rightarrow(1)}$: Assume that the series (\ref{Eq:W:1:12}) diverges and prove that 
\begin{equation}\label{Eq:W:3:12}
H_{1_\infty}^{\bf \Omega} \equiv 0
\end{equation} 
Proof is based on the construction of the family
$\{G_p\}$ of nonnegative ${\cal A}$-harmonic functions in $\Omega$ with the following properties: 
\begin{equation}\label{Eq:W:3:13}
\lim_{x\to \infty} G_p(x)=1 \quad\text{ for any fixed}\ p
\end{equation}
\begin{equation}\label{Eq:W:3:14}
\lim_{p\to +\infty} G_p(x)=0 \quad\text{ for any fixed}\ x \in \Omega
\end{equation}
Indeed, in this case family $\{G_p\}$ is in the upper class ${\cal S}(\Omega)$
with respect to the Perron solution $H_{1_\infty}^{\bf \Omega}$. Accordingly, we have
\[ 0\le H_{1_\infty}^{\bf \Omega}(x) \le G_p(x), \ x\in \Omega  \]
and passing to the limit as $p\to +\infty$, (\ref{Eq:W:3:12}) follows.

To construct the required family of ${\cal A}$-harmonic functions, we need to enter
the positive integer parameter $p$ into the Wiener series and make some rearrangements. Let
\[ \Gamma_n^p\equiv cap(E_n^p), E_n^p \equiv \Omega^c \cap \{2^{\frac{n-1}{p}} \le |x| \le 2^{\frac{n}{p}} \}  \]
Assuming that $n\ge p$, from the subadditivity of the capacitary measure it follows that
\[ \Gamma_n \le \Gamma_{np}^p+\Gamma_{np-1}^p+\cdots+\Gamma_{np-n+1}^{p}  \]
Accordingly, one of the following series must be divergent:
\[ \sum_n 2^{-n(N-2)}\Gamma_{np}^p \ \ \sum_n 2^{-n(N-2)}\Gamma_{np-1}^p \ \cdots \   \sum_n 2^{-n(N-2)}\Gamma_{np-n+1}^p  \]
It is easy to see that the divergence of any of these series 
implies that
\[ \sum_n 2^{-\frac{n(N-2)}{p}}\Gamma_n^p=+\infty.  \]
By choosing $p^2$ successive values of $n$  as
\[ kp^2, kp^2+1,...,kp^2+p^2-1, \]
it follows that one of the following series must also be divergent:
\[ \sum_n 2^{-np(N-2)}\Gamma_{np^2}^p \ \ \sum_n 2^{-(np+\frac{1}{p})(N-2)}\Gamma_{np^2+1}^p \ \cdots \   \sum_n 2^{-((n+1)p-\frac{1}{p})(N-2)}\Gamma_{(n+1)p^2-1}^p  \]
We may therefore assume without loss of generality, that
\begin{equation}\label{Eq:W:3:15}
 \sum_n 2^{-np(N-2)}\Gamma_{np^2}^p=+\infty.
\end{equation}
The proof given for the contrary case is similar to the one presented. 
In view of (\ref{Eq:W:3:15}), for an arbitrary large integer $p$ we can find 
an integer $N_p$ such that
\begin{equation}\label{Eq:W:3:16}
 \sum_{n=1}^{N_p} 2^{-np(N-2)}\Gamma_{np^2}^p>p.
\end{equation}
If $p>2$, the distance between $E_{np^2}^p$ and its closest neigbors may be estimated as follows
\[ dist (E_{np^2}^p; E_{(n-1)p^2}^p) \ge 2^{np}(2^{-\frac{1}{p}}-2^{-p}), \]
\[ dist (E_{np^2}^p; E_{(n+1)p^2}^p) \ge 2^{np}(1-2^{-\frac{p}{2}}). \]
Accordingly, we have
\[ dist (E_{np^2}^p; E_{kp^2}^p)\ge 2^{np} \alpha(p), \quad\text{ for }\ \forall k \neq n,  \]
where
\[ \alpha(p)=min(2^{-\frac{1}{p}}-2^{-p}; 1-2^{-\frac{p}{2}}).  \]
By using (\ref{Eq:W:2:10}) we have
\begin{equation}\label{Eq:W:3:17}
 \zeta_{E_{np^2}^p}(x)\le K 2^{-np(N-2)}\Gamma_{np^2}^p \alpha^{2-N}(p), \ x\in E_{kp^2}^p, k\neq n. 
\end{equation}
Since the function $\sum_{n=1}^{N_p} \zeta_{E_{np^2}^p}$ is ${\cal A}$-harmonic in $\Big ( \cup_{n=1}^{N_p}E_{np^2}^p \Big )^c$ including in $\Omega$, and vanishes at $\infty$,
from (\ref{Eq:W:3:17}) it follows that
\begin{equation}\label{Eq:W:3:18}
\sum_{n=1}^{N_p} \zeta_{E_{np^2}^p}(x)\le K + K \alpha^{2-N}(p) \sum_{n=1}^{N_p}2^{-np(N-2)}\Gamma_{np^2}^p, \ x\in \Big ( \cup_{n=1}^{N_p}E_{np^2}^p \Big )^c . 
\end{equation}
Now, we choose the required family of ${\cal A}$-harmonic functions as follows:
\[ G_p(x)=-\frac{\sum_{n=1}^{N_p} \zeta_{E_{np^2}^p}(x)}{\gamma_p \sum_{n=1}^{N_p}2^{-np(N-2)}\Gamma_{np^2}^p} + 1, \ \gamma_p = K(p^{-1}+\alpha^{2-N}(p)). \]

Nonnegativity of $G_p$ follows from (\ref{Eq:W:3:16}) and (\ref{Eq:W:3:18}):
\[ G_p(x)\ge \frac{K\Big(-1-\alpha^{2-N}(p)\sum_{n=1}^{N_p}2^{-np(N-2)}\Gamma_{np^2}^p \Big )}
{\gamma_p \sum_{n=1}^{N_p}2^{-np(N-2)}\Gamma_{np^2}^p } + 1 \]

\[ \ge -\frac{K(p^{-1}+\alpha^{2-N}(p))}
{\gamma_p}+1=0, x\in \Omega. \]

Since 
\[ \lim_{x\to \infty} \zeta_{E_{np^2}^p(x)}=0, \ n=1,2,...,N_p, \]
(\ref{Eq:W:3:13}) easily follows. \\

Let $x\in \Omega$ is fixed. From (\ref{Eq:W:2:10}) it follows that
\[ \zeta_{E_{np^2}^p}(x)\ge \frac {K^{-1} \Gamma_{np^2}^p}{|2^{np}+|x||^{N-2}}
\ge \frac{K^{-1}2^{-np(N-2}\Gamma_{np^2}^p}{|1+2^{-p}|x||^{N-2}}. \]  

Hence we have
\[ G_p(x)\le -\frac{K^{-1}}{\gamma_p |1+2^{-p}|x||^{N-2}}+1, \]
and (\ref{Eq:W:3:14}) immediately follows. Regularity of $\infty$ is proved.\\

${\bf (1)\Rightarrow(5)}$: Assume that $\Omega^c$ is ${\cal A}$-thin at $\infty$.
Prove that $\infty$ is irregular for $\Omega$. We will only consider the
non-trivial case when $\Omega^c$ is unbounded. It is sufficient to show that
$\infty$ is irregular for $\Omega \cup \Sigma$, where $\Sigma$ is some open sphere
of large radius. This immediately follows from the equivalence of the 
irregularity of $\infty$ to convergence of the Wiener series (\ref{Eq:W:1:12}).
Obviously the latter is not affected by adding a finite number of terms.

By Definition 1.2 and (\ref{Eq:W:1:11}), it is easy to construct an 
${\cal A}$-superharmonic function $u$ such that
\[ 0\equiv u(\infty)<2\equiv\liminf_{x\to\infty, x\in \Omega^c}u(x) \]
Choose an open sphere $\Sigma$ of large enough radius such that
\[ u(x)\ge 1 \quad\text{ for }\ x\in \Omega^c \cap \Sigma^c \]
It is easy to see that ${\cal A}$-subharmonic function $v=1-u$
satisfies
\[ \limsup_{x\to \infty} v =1  \]
and belongs to lower class $-{\cal S}(\Omega\cup \Sigma)$ with
respect to the Perron solution $H_{1_\infty}^{\bf \Omega \cup \Sigma}$.
Therefore we have
\[ v(x)\le  H_{1_\infty}^{\bf \Omega \cup \Sigma}(x) \quad\text{ for }\ x\in \Omega \cup \Sigma \]
which implies that
\[  \limsup_{x\to \infty}H_{1_\infty}^{\bf \Omega \cup \Sigma}(x)  =1. \]
Hence $\infty$ is irregular for $\Omega \cup \Sigma$, and for $\Omega$ as well.\\

${\bf (5)\Rightarrow(1)}$: Assume that $\infty$ is irregular for $\Omega$. Prove
that $\Omega^c$ is ${\cal A}$-thin at $\infty$. As before, we will only consider the
non-trivial case when $\Omega^c$ is unbounded. It is sufficient to demonstrate
that $(\Omega \cup \Sigma)^c$ is ${\cal A}$-thin at $\infty$, where $\Sigma$ is a sphere 
of large radius. 

While proving the relation ${\bf (1)\Rightarrow(4)}$ above, we constructed function
$\vartheta$ which is ${\cal A}$-superharmonic in $\rr^N$ and satisfies (\ref{Eq:W:3:6}),
and 
\[ \vartheta \ge 1 \quad\text{ quasieverywhere on }\ (\Omega^m)^c\equiv \Big ( 
\bigcup_{n=m+1}^{+\infty} E_n \Big )   \]
We want to have this property everywhere on $(\Omega^m)^c$. 
Let
\[ e_n = \{x\in E_n\cap \Sigma_{2^n}: \vartheta(x)<1 \}, \ n=m+1,m=2,\cdots   \]
Since $cap_{\cal A}(e_n)=0$, for arbitrary positive $\delta_n$, we can choose
open cover $c_n$ of $e_n$ such that $cap_{\cal A}(d_n)<\delta_n$, where $d_n=\ov{c_n}$. 
We can obviously choose open sets $c_n$ in such a way that 
\[ c_n\subset \{x: |x|\ge 2^{n-\frac{1}{2}} \}. \]
Clearly, $x^*$
lies outside $\Big ( 
\bigcup_{n=m+1}^{+\infty} E_n \cup c_n \Big )$. Consider
${\cal A}$-superharmonic function
\[ \varsigma(x) \equiv \sum_{n=m+1}^{+\infty}\zeta_{d_n}(x) \]
We estimate $\varsigma(x^*)$ similar to the estimation of $\vartheta(x^*)$.
From (\ref{Eq:W:2:8})-(\ref{Eq:W:2:10})
it follows that
\[ \zeta_{d_n}(x^*)\le K \int_{d_n} |y-x^*|^{2-N}d\lambda_n(y) \]
\[ \le K \Big ( \frac{4}{\sqrt{2}-1} \Big )^{N-2} 2^{-n(N-2)}cap_{\cal A}(d_n), 
\quad\text{ for }\ n\ge m+1, \]
\[ \varsigma(x^*)\le K \Big ( \frac{4}{\sqrt{2}-1} \Big )^{N-2} \sum_{n=m+1}^{+\infty} 2^{-n(N-2)}\delta_n \]
where $\lambda_n$ is an ${\cal A}$-capacitary measure of $d_n$. By choosing
numbers $\delta_n$ sufficiently small, we have
\[ \varsigma(x^*) < \frac{\ep}{2} \]
Consider a function $\varpi=\varsigma+\vartheta$. It is ${\cal A}$-superharmonic 
in $\rr^N$,
\[ \varpi (x)\ge 1 \quad\text{ everywhere on }\ \ (\Omega^m)^c. \]
and 
\[ \varpi (x^*)<\ep  \]
Since $\varpi$ is ${\cal A}$-superharmonic in $\Omega^m$, by applying the
${\cal A}$-superharmonic minimum principle in $\Omega^m \cap \{x: |x|<R\}$ we conclude 
that for all sufficiently large $R$, there must be a point $x_R\in \Omega_e\cap\{x:|x|=R\}$ such that
\[ \varpi(x_R)<\ep \] 
Hence, we deduce that
\[ \liminf_{x\to \infty} \varpi(x) \le \ep. \]
Therefore, $\varpi$ is an ${\cal A}$-superharmonic function which guarantees
the ${\cal A}$-thinness of $\Omega^c$ according to the Definition 1.2 and 
(\ref{Eq:W:1:11}).



\end{document}